\documentclass[oneside]{amsart}
\usepackage{amsfonts}
\usepackage{amssymb}
\usepackage{amsxtra}
\usepackage{amstext}
\usepackage[english]{babel}

\newcommand{\R}{\mathbb{R}}
\newcommand{\C}{\mathbb{C}}
\newcommand{\D}{\mathbb{D}}
\newcommand{\N}{\mathbb{N}}
\newcommand{\E}{\mathbb{E}}
\newcommand{\Z}{\mathbb{Z}}

\newcommand{\pp}{\mathbb{P}}
\newcommand{\kA}{\mathcal{A}}

\newcommand{\kR}{\mathcal{R}}

\newcommand{\kP}{\mathcal{P}}

\newcommand{\kF}{\mathcal{F}}
\newcommand{\kE}{\mathcal{E}}

\newcommand{\kN}{\mathcal{N}}
\newcommand{\kX}{\mathcal{X}}

\newcommand{\ka}{\mathfrak{a}}

\newtheorem {lem} {Lemma} [section]
\newtheorem {prop} {Proposition} [section]
\newtheorem {theo} {Theorem} [section]

\newtheorem {rem} {Remark} [section]

\newcommand{\la}{\lambda}
\newcommand{\teta}{\theta}

\title[Random walk on a building of type $\tilde{A}_r$ and Brownian motion of the Weyl chamber]
      {Random walk on a building of type $\tilde{A}_r$ and Brownian motion of the Weyl chamber}
\author{Bruno Schapira}

\begin{document}

\begin{abstract} In this paper we study a random walk on an affine building of type $\tilde{A}_r$, whose
radial part, when suitably normalized, converges toward the Brownian
motion of the Weyl chamber. This gives a new discrete approximation
of this process, alternative to the one of Biane \cite{Bi2}. This
extends also the link at the probabilistic level between Riemannian
symmetric spaces of the noncompact type and their discrete
counterpart, which had been previously discovered by Bougerol and
Jeulin in rank one \cite{BJ1}. The main ingredients of the proof are
a combinatorial formula on the building and the estimate of the
transition density proved in \cite{AST}.
\end{abstract}

\maketitle
\let\languagename\relax


\bigskip

\noindent \textbf{A.M.S. classification:} 05C25, 60B10, 60B15, 60C05,
60J10, 60J25, 60J35, 60J60.

\section{Introduction}

The Brownian Motion of the Weyl chamber, or intrinsic BM, considered by Biane in \cite{Bi1}, has bring recently a growing attention, due to his apparition in many branches of the probability theory. In dimension $1$ for instance it is the Bessel 3 process, whose importance is well known.  
More generally it appears in particular in the theory of random matrices, 
of particles systems, in queuing theory, in oriented percolation (see \cite{O} for an overview), but also in the theory of stochastic processes on Riemannian symmetric spaces \cite{ABJ}, and more recently in the theory of Dunkl and Heckman-Opdam processes \cite{GY} \cite{S}. It takes values in a Weyl chamber, denoted here by $\ka_+$, which is a cone of the Euclidean space $\R^r$ delimited by hyperplanes satisfying some conditions. More precisely it may be defined as the Brownian motion killed on the boundary of $\ka_+$ and conditioned (in the sense of Doob) never to escape the cone. 

\vspace{0.2cm}

A basic question is to find a natural discrete version of this process, i.e. some random walk on some lattice of $\R^r$ which after renormalization converges toward this process. An example of such random walk has been introduced and first studied by Biane \cite{Bi2} \cite{Bi3} \cite{Bi4} and then by Biane, Bougerol and O'Connell \cite{BiBO}. It is defined on the so-called weight lattice, actually its positive part $P^+$, which is a natural choice, well known in the theory of root systems, since it has roughly the same geometric structure than $\ka_+$ (see the next section for details). So this random walk is defined analogously, as the simple random walk on $P^+$, killed on the boundary of $P^+$, and conditioned never to touch it. One can compute explicitly its transition kernel, with the help of the reflexion principle, and then prove its convergence, after the usual renormalization, toward the intrinsic BM \cite{BiBO}. In fact \cite{Bi2} \cite{Bi4} this random walk appears naturally in the context of quantum random walks on the dual of a compact Lie group.

\vspace{0.2cm}

In this paper we study another discrete version of the intrinsic BM, when the cone $\ka_+$ is associated to a root system of type $A_r$ (see the next section). In this case the intrinsic BM identifies with the process of eigenvalues of hermitian matrices with zero trace, whose coefficients are Brownian motions (also called GUE process). Our random walk is also defined on the weight lattice but indirectly. It is the so-called radial part of some nearest neighbor random walk $(X_n,n\ge 0)$ on an affine building of type $\tilde{A}_r$. To be short, let say simply that a building is some graph containing copies, called apartments, of the weight lattice. Once a vertex $O$, an apartment containing $O$, and its positive part $P^+$, have been fixed, then the radial part of a vertex $x$ is its projection $\overline{x}$, in some sense, onto $P^+$. This notion of radial part generalizes the usual notion of distance between two vertices in regular trees, which are buildings of type $\tilde{A}_1$. Call $A$ the transition operator of the random walk $(X_n,n\ge 0)$. Let $F_0$ be some eigenfunction of $A$ at the bottom of its spectrum. Let $(Y_n,n\ge 0)$ be the relativized $F_0$-random walk in the sense of Doob. For $N\ge 0$, let $(Y^N_t,t\ge 0)$ be the continuous time normalized $F_0$-random walk defined by $Y_t^N=\overline{Y_{[Nt]}}/\sqrt{N}$ for all $t\ge 0$. Then the main result of this paper is that the sequence $(Y^N_t,t\ge 0)$ converges in law in the path space toward the intrinsic BM, when $N\to +\infty$.

The idea of this result goes back to the work of Bougerol and Jeulin \cite{BJ2}. They proved it in the tree case, which as we already mentioned are particular examples of buildings (they are affine buildings of rank one). In the same time they proved an analogue version on Riemannian symmetric spaces of the noncompact type of rank one (in particular on complex and real hyperbolic spaces), which are continuous versions of regular trees. Then they managed to extend their result with Anker \cite{ABJ} \cite{BJ1} on Riemannian symmetric spaces of higher rank. The author also proved it in the general context of Heckman--Opdam's theory \cite{S}. So the present paper investigates the higher rank case on the discrete space level. However this story may not be finished yet, since here we consider essentially only one particular random walk and only on a particular type of affine buildings (those associated to root systems of type $A$). The main reason for this is that we need good estimates of the transition densities of the random walk, which for the moment are known only in few cases \cite{AST}. Random walks on buildings have already been studied so far. In particular \cite{CaW} and \cite{Pa3} prove a law of large numbers, central and local limit theorems.

\vspace{0.2cm}

This paper is organized as follows. In the next section we recall some definitions and basic properties of the main objects of study. In section \ref{sectionF0} we define precisely the $F_0$-random walks and show that they appear naturally as limits of random bridges when the length of the bridge tends to infinity. Then in section \ref{seccombi}, we compute explicitly the transition densities of the radial part of nearest neighbor random walks. This, in our opinion, is actually the main part of the proof (once the results of \cite{AST} are known). It allows to derive directly the result for all nearest neighbor random walks when the starting point of the intrinsic BM lies inside the Weyl chamber. For this we have to renormalize also the starting point of the random walks and use a general criteria of Ethier and Kurtz \cite{EK} on the convergence of generators. We explain all this in section \ref{secdepartbiaise}. In the last section we consider the case where the intrinsic BM starts from the origin. This case is more difficult, since the criteria of Ethier and Kurtz is not sufficient. We need in addition estimates of the transition densities of the random walk to control its norm in small times. These estimates, as we already mentioned were proved recently in \cite{AST}.

\section{Preliminaries}
\noindent \textbf{The root system.} We denote by $\ka$ the Euclidean space of dimension $r$, endowed with its scalar product $<\cdot,\cdot>$. We denote by $\ka_\C:=\ka+i\ka$ its complexification. Let $\kR\subset \ka$ be a root system of type $A_r$. One can define it (see \cite{Bou}) as the set of vectors $e_i-e_j$, for $1\le i\neq j \le r+1$, where $(e_1,\dots, e_{r+1})$ is the canonical basis of $\R^{r+1}$ and $\ka$ is the hyperplane $\{<x,e_1+\dots +e_{r+1}>=0\}$, endowed with the restriction of the canonical scalar product on $\R^{r+1}$. We denote by $|\cdot|$ the associated norm. A subset of positive roots $\kR^+$ is defined as follows: one first choose arbitrarily $u\in \ka$ such that $<u,\alpha>\neq 0$ for all $\alpha \in \kR$. Then one set $\kR^+=\{\alpha \in \kR \mid <\alpha,u>\ > 0\}$. For instance one can take the set of vectors $e_i-e_j$, where $i<j$. The positive Weyl chamber associated to $\ka_+$ is defined by 
$$\ka_+=\{x \in \ka \mid <\alpha,x> \ >0 \quad \forall \alpha \in \kR^+\}.$$
We denote by $\overline{\ka_+}$ its closure and by $\partial\ka_+$ its boundary. If $\alpha \in \kR$, we denote by $r_\alpha$ the orthogonal reflexion with respect to the hyperplane orthogonal to $\alpha$. The Weyl group $W_0$ is the finite subset of the orthogonal group generated by the $r_\alpha$'s with $\alpha\in  \kR$. In fact $\kR$ is generated by the set $\Pi:=\{\alpha_1,\dots,\alpha_r\}$, where $\alpha_i=e_i-e_{i+1}$ (these elements are called the simple roots). Furthermore $W_0$ is generated by the $r_i:=r_{\alpha_i}$, $i=1,\dots,r$. Any $w\in W_0$ can be written under the form $r_{i_1}\dots r_{i_l}$. The smallest possible $l$ in such writing is called the length of $w$, and is denoted by $l(w)$. The dual basis of $\Pi$ is denoted by $\{\la_1,\dots,\la_r\}$. More precisely we assume that $<\alpha_i,\la_j>=\delta_{i,j}$ for all $i,j$. The elements of the dual basis are called the fundamental weights. The weight lattice $P$ is by definition the $\Z$-lattice generated by the fundamental weights. We denote by $P^+$ the subset of positive weights, i.e. $P\cap \overline{\ka_+}$, and by $P^{++}$ the subset of regular positive weights, i.e. $P\cap \ka_+$. The lattice $P$ generates a simplicial complex, called the Coxeter complex (see \cite{Pa1} for a more precise definition), which is denoted by $\kP$. We attribute labels lying in the set $\{0,\dots,r\}$ to the elements of $P$ as follows: we attribute the label $i$ to $\la_i$, and for any simplex of maximal dimension we assume that the set of labels of its extremal vertices is exactly $\{0,\dots,r\}$. We denote by $W$ the affine Weyl group, which is generated by $W_0$ and the translations by elements of $P$. We denote by $\widetilde{W}$ the extended affine Weyl group, which is generated by $W$ and the group $\Omega$ which acts by circular permutation on the fundamental weights and $0$ (see \cite{Bou}). The group $W$ is in fact generated by $\Pi$ and the orthogonal reflexion $r_0$ with respect to the hyperplane $\{<x,\alpha_1+\dots +\alpha_r>\}=1$. This allows to define the length of $w\in W$ analogously as in $W_0$. Then we set $q_w=q^{l(w)}$. This definition extends to $\widetilde{W}$: if $\tilde{w}=wg$, with $w\in W$ and $g\in \Omega$, then we set $q_{\tilde{w}}=q_w$. One can prove (see \cite{Pa1}) that if $t_\la$ is the translation by $\la\in P^+$, then 
$$q_{t_\la}=q^{\sum_{\alpha\in \kR^+}<\alpha,\la>}.$$
It will be convenient to extend this formula to any $\la \in P$. So we set: 
$$
\tilde{q}_{t_\la}:=q^{\sum_{\alpha\in \kR^+}<\alpha,\la>} \quad \forall \la \in P.
$$

\noindent \textbf{The affine building.} An affine building (see \cite{Pa1} or \cite{Ron}) of type $\tilde{A}_r$ is a nonempty simplicial complex containing sub-complexes, called apartments, such that: 
\begin{itemize}
\item Any apartment is isomorphic (see \cite{Pa1}) to the Coxeter complex $\kP$. 
\item Given two chambers (simplex of maximal dimension) there exists an apartment containing both. 
\item Given two apartments having at least one common chamber there exists a unique isomorphism between them fixing pointwise their intersection.
\end{itemize}
We denote by $\kX$ the set of vertices (simplexes of dimension $0$) of the building and we fix one, called $O$. The preceding definition allows to attribute labels to elements of $\kX$ such that $O$ has label $0$ and isomorphisms in the definition preserve labels.

Let us mention, to avoid any ambiguity, that chambers of the building, also called alcoves, which are compact simplexes, must not be confused with the positive Weyl chamber, which is a chamber of the root system and in particular an unbounded open cone of $\ka$.

Two chambers of the building are called adjacent if they are distinct and have a common face (simplex of co-dimension $1$). A gallery is a sequence of chambers $(C_0,\dots,C_n)$, where for all $i\le n-1$, $C_i$ and $C_{i+1}$ are adjacent. The integer $n$ is called the length of the gallery. We define the distance between two chambers $C$ and $C'$ as the minimal length of a gallery starting from $C$ and ending in $C'$, and we denote it by $d(C,C')$. If $F$ a face of a chamber $C$, we will assume that the number $q$ of chambers adjacent to $C$ and containing $F$ is independent of $F$ and $C$ (we say that the building is regular). For $r=1$ for instance, $\kX$ is a tree where each vertex has $q+1$ neighbors. 

\vspace{0.2cm}

If $x\in \kX$ there exists by definition an apartment containing $x$ and $O$, and an isomorphism between this apartment and $\kP$ sending $O$ on $0$ and $x$ on an element of $P^+$. The image of $x$ by this isomorphism, let say $\overline{x}$, is uniquely defined and called the radial part or coordinate of $x$. When $r=1$ for instance $P^+$ identifies with $\N$ and $\overline{x}$ is the distance (in the usual sense) between $x$ and $O$. For $\la \in P^+$, we denote by $V_\la(O)$ the set of vertices of coordinate $\la$ and we call it the sphere of radius $\la$ centered in $O$. Given a vertex $x$ we define also the sphere $V_\la(x)$ of radius $\la$ centered in $x$, as the set of vertices $y\in \kX$ such that there exists an apartment $\kA$ containing $x$ and $y$ and an isomorphism between $\kA$ and $\kP$, preserving labels up to translation, and sending $x$ on $0$ and $y$ on $\la$. For all $\la \in P^+$, $V_\la(x)$ is a finite subset of $\kX$, whose cardinality $N_\la$ is independent of $x$. If $\la \in P$, let $W_{0\la}$ be the stabilizer of $\la$ under the action of $W_0$. Then we have the formula (see \cite{Pa1} Formula (1.5)): 
$$N_\la = \frac{W_0(q^{-1})}{W_{0\la}(q^{-1})}q_{t_\la},$$
where $V(q^{-1})=\sum_{w\in V} q_w^{-1}$ for any subgroup $V$ of $W_0$. We define also the function $\pi$ on $P$ by 
$$\pi(\la)= \prod_{\alpha \in \kR^+} <\alpha,\la>.$$

\noindent \textbf{The Macdonald polynomials.} The functions
$\mathbf{c}$ and $h$ are defined for $z\in \ka_\C$ by
$$\textbf{c}(z)=\prod_{\alpha\in
\kR^+}\frac{1-q^{-1}e^{-<\alpha,z>}}{1-e^{-<\alpha,z>}},$$
and
$$h(z)=\sum_{i=1}^r\sum_{\la\in W_0\la_i}e^{<\la,z>}.$$ 
We set also $\tilde{h}:=h/h(0)$. Macdonald's polynomials are defined (see for instance \cite{Ca} \cite{M1} \cite{M2} or
\cite{Pa1}) for $\la \in P^+$ and $z\in \ka_\C$ by
\begin{eqnarray}
\label{polynomemac}
P_\la(z)=\frac{q_{t_\la}^{-{\frac{1}{2}}}}{W_0(q^{-1})}\sum_{w\in
W_0}\mathbf{c}(w^{-1} z) e^{<w\la,z>}.
\end{eqnarray}
The function $F_0$ is defined on $P^+$ by
$$F_0(\la)=P_\la(0).$$
We recall the estimate in $P^+$ (see \cite{AST}):
\begin{eqnarray}
\label{estimationF01} F_0(\la)\asymp
q_{t_\la}^{-\frac{1}{2}}\prod_{\alpha \in \kR^+}(1+<\alpha,\la>).
\end{eqnarray}
Moreover when $<\alpha,\la>\to +\infty$ for all $\alpha \in \kR^+$, then
\begin{eqnarray}
\label{estimationF02} F_0(\la)\sim
\textrm{const}\cdot q_{t_\la}^{-\frac{1}{2}}\pi(\la).
\end{eqnarray}

\noindent \textbf{Symmetric nearest neighbor random walk.} By definition it is a Markov chain $(X_n,n\ge 0)$ on $\kX$ whose transition densities are equal to
\begin{eqnarray*}
p(x,y)=\left\{ \begin{array}{cl}
 p_i & \text{if }y\in V_{\la_i}(x),\\
0 & \text{otherwise},
\end{array}
\right.
\end{eqnarray*}
where the $p_i$'s satisfy the condition $\sum_{i=1}^r p_iN_{\la_i}=1$.
Let
$$\tilde{\rho}=\sum_{i=1}^rp_iq_{t_{\la_i}}^{\frac{1}{2}}|W_0\la_i|,$$
be the spectral gap (the fact that $\tilde{\rho}$ is well the spectral gap results from Formula (2) in \cite{AST}). 
The radial random walk is by definition the Markov chain $(\overline{X_n},n\ge 0)$ on $P^+$. 
We denote by $\overline{p}$ its transition kernel. A particular example is the so-called simple random walk. Its transition probabilities are determined by
\begin{eqnarray}
\label{probatransition}
p_i=\frac{q_{t_{\la_i}}^{-\frac{1}{2}}}{\sum_{i=1}^rq_{t_{\la_i}}^{-\frac{1}{2}}N_{\la_i}},
\end{eqnarray}
for all $i\le r$. For this random walk we have a relatively simple integral formula for the transition kernel (see for instance \cite{Pa3} Formula (1.10) and \cite{AST} Formula (2)): 
if
$$U=\{\teta \in \ka \mid <\alpha,\teta>\ \le
2\pi \quad \forall \alpha \in \kR \},$$ and if $p_n(O,x)$ denotes the
probability, starting from $O$, to arrive in $x\in V_\la(O)$
in $n$ steps, then
\begin{eqnarray}
\label{formuleintegtransition} p_n(O,x)=\text{const}\cdot
\tilde{\rho}^n
\int_U\tilde{h}^n(i\teta)P_\la(i\teta)\ \frac{d\teta}{|\mathbf{c}(i\teta)|^2}.
\end{eqnarray}

\noindent \textbf{The Brownian motion of the Weyl chamber.} This process is also called the intrinsic Brownian motion. There are many ways to define it. We will only recall some of them here.   
Before this, to avoid some (linear) time change in Theorem \ref{theoprincipal}, we define another norm on $\ka$:   
$$||x||^2:=\frac{1}{h(0)}\sum_{i=1}^r\sum_{\la \in W_0 \la_i} <\la,x>^2.$$
Since any homogeneous $W_0$-invariant polynomial of degree $2$ is proportional to $|\cdot|^2$ (see \cite{Hel}), one knows that $||\cdot||$ is proportional to $|\cdot|$. An elementary calculus shows that in fact $||\cdot||^2=2^{r-2}(2^r-1)^{-1}|\cdot|^2$. We will denote by $<<\cdot,\cdot>>$ the scalar product associated to $||\cdot||$, and we set $c:=2^{r-2}(2^{r}-1)^{-1}$.

Now one can define the intrinsic BM (see Biane \cite{Bi1} or Anker--Bougerol--Jeulin \cite{ABJ}) as the $\pi$-process, in the sense of Doob, of the standard Brownian motion in $(\ka,||\cdot||)$ killed on the boundary of $\ka_+$. In fact, as Biane noticed, $\pi$ is the unique positive harmonic function which vanishes on $\partial \ka_+$. So one can also interpret this process as the Brownian motion killed on $\partial \ka_+$ and conditioned to never touch $\partial \ka_+$ (or to escape $\ka_+$ at infinity). From this point of view it is a natural generalization of the Bessel-3 in higher dimension. In fact Biane, Bougerol and O'Connell \cite{BiBO} proved a deeper result which reinforces this. They proved that this process can be obtained by some transformation of the Brownian motion on $\ka$, which generalizes the Pitman transform $2S-B$ in dimension $1$. But one can also define it as the Dunkl process with parameter $k=1$ (see \cite{GY}). In particular it is a Feller process with generator $D$ (see \cite{S}) defined for a regular and $W_0$-invariant function $f$ by   
$$Df(x)=\frac{1}{2}\Delta f(x) + <<\nabla \log \pi, \nabla f>>(x) \quad \forall x \in \ka_+,$$
where $\Delta$ is the usual Laplacian on $(\ka,||\cdot||)$. Finally, since our definition coincides with the definition in \cite{ABJ} up to a linear time change, one can deduce from Formula (3.5) in 
\cite{ABJ} that the law of this process starting from $0$ has density
$$p_t(0,x)=\text{const}\cdot
\frac{1}{t^{|\kR^+|+\frac{r}{2}}}\pi(x)^2 e^{-\frac{|x|^2}{2ct}},$$
for all $t\ge 0$ and all $x\in \ka_+$.

\section{The $F_0$-random walk}
\label{sectionF0}
Let $(X_n,n\ge 0)$ be a symmetric nearest neighbor random walk, with transition kernel $p$. The function $F_0$ is an eigenfunction of its transition operator (see \cite{Pa1} Theorem 3.22):
$$\sum_{y\in \kX}p(x,y)F_0(\overline{y})=\tilde{\rho}F_0(\overline{x}),$$
for all $x\in \kX$. Then the $F_0$-random walk 
$(Y_n,n\ge 0)$ is the Markov chain on $\kX$ with transition kernel $q$ 
$$q(x,y)=p(x,y)\frac{F_0(\overline{y})}{F_0(\overline{x})}\tilde{\rho}^{-1}.$$
We denote by $(X_n^{N,0},n\ge 0)$ the bridge of length $N$ around $O$, i.e. the random walk starting from $O$ and conditioned to come back in $O$ at time $N$. 
The next proposition shows that the $F_0$-random walk can be seen as a loop of infinite length around $O$. There is an analogue result on Riemannian symmetric spaces $G/K$ \cite{ABJ} and in Heckman--Opdam's theory \cite{S}.  

\begin{prop} \label{convergencepont} When $N\to +\infty$, $(X_n^{N,0},n\ge 0)$ converges
in law, in the path space, toward $(Y_n,n\ge 0)$.
\end{prop}
\begin{proof} Let $\pp$ be the law of $(X_n,n\ge 0)$ and let 
$\kF_n:=\sigma(X_k,k\le n)$ be its natural filtration. Let
$\pp^{N,0}$ be the law of $(X_n^{N,0},n\ge 0)$. The following absolute continuity relation holds:
\begin{eqnarray}
\label{abscontinuite} \pp^{N,0}_{|\kF_n}=
\frac{p_{N-n}(O,X_n)}{p_N(O,O)}\pp_{|\kF_n},
\end{eqnarray}
for all $n\ge 0$. To simplify, assume first that $(X_n,n\ge0)$ is the simple random walk. Then from 
\eqref{formuleintegtransition}, we deduce that for all $\la\in P^+$ and all
$x\in V_\la(O)$,
\begin{eqnarray*}
\frac{p_{N-n}(O,x)}{p_N(O,O)}=\tilde{\rho}^{-n}\frac{\int_U
\tilde{h}^{N-n}(i\teta)P_\la(i\teta)|\mathbf{c}(i\teta)|^{-2}\ d\teta}{\int_U
\tilde{h}^{N}(i\teta)P_0(i\teta)|\mathbf{c}(i\teta)|^{-2}\ d\teta}.
\end{eqnarray*}
Next in the above integrals we make the following change of variables: $\teta\to \teta/\sqrt{N-n}$, respectively
$\teta \to \teta/\sqrt{N}$. Then one can observe that 
$$\tilde{h}^N(i\frac{\teta}{\sqrt{N}})\to
e^{-\frac{1}{2 h(0)} \sum_\la <\la,\teta>^2}=e^{-\frac{1}{2}||\teta||^2} ,$$ when $N\to
+\infty$. Moreover $\mathbf{c}^{-1}(i\teta)$ is equivalent (up to some constant) to $\pi(i\teta)$ near $0$. Thus we get the following convergence result: 
$$\frac{p_{N-n}(O,x)}{p_N(O,O)}\to
\tilde{\rho}^{-n}\frac{P_\la(0)}{P_0(0)},$$ when $N\to +\infty$. Together with
the relation \eqref{abscontinuite} this proves the proposition for the simple random walk. In general $\tilde{h}$ has to be replaced by some sum of exponentials but whose coefficients are not necessarily equal (see \cite{Pa3} Formula (1.10) and \cite{AST} Formula (2)). However the same proof applies as well. 
\end{proof}

\section{Transition probabilities of the radial random walk}
\label{seccombi}
In this section we give explicit formula for the transition probabilities of radial part of nearest neighbor random walks. We will use them in the next section to prove Theorem \ref{theogeneral}. These formula are derived from a combinatorial calculus (Proposition \ref{lemprobatransition} below) which could be of independent interest. One can find a similar result in \cite{CaW} Lemma 2.1, and at the end of Parkinson's thesis \cite{Pa2} for all affine buildings of rank $2$.  

\begin{prop}
\label{lemprobatransition} For all $\la \in P^{++}$, 
$x_\la \in V_\la(O)$, $w\in W_0$ and $i\in
\{1,\dots,r\}$,
$$\left|V_{\la+w\la_i}(O)\cap
V_{\la_i}(x_\la)\right|=q_{t_{\la_i}}^{\frac{1}{2}}\tilde{q}_{t_{w\la_i}}^{\frac{1}{2}}.$$
\end{prop}
\begin{proof} Let us fix $\la \in P^{++}$ and $x_\la\in V_\la(O)$.
Let $w'\in W_0$ be of maximal length such that $w'\la_i=w\la_i$. Since the desired formula does not change when we replace $w$ by $w'$, one can always assume that $w=w'$. 
Now by definition of $q$, we know that if we fix a face $F$ containing $x_\la$, then there are exactly $q+1$ chambers containing $F$. Equivalently, if we fix a chamber $C$ and one of its faces $F$, there are exactly $q$ vertices $x\in \kX$, such that $(F,x)$ defines a chamber adjacent to $C$. Denote by $C_+$ and $C_-$ the chambers (in $\kP$) containing $\la$ and all weights $\la+\la_i$, respectively $\la-\la_i$, for $i=1,\dots,r$. We set also $C_w:=t_\la w t_{-\la}C_+$. We know (see \cite{Pa1} Lemma B.2) that for all $i$, 
\begin{eqnarray} \label{convexite}
\left|V_{\la-\la_i}(O)\cap V_{\la_j}(x_\la)\right|=1, \end{eqnarray}
if $-\la_i\in W_0\la_j$. In other words there exists a unique chamber in the building, also denoted by $C_-$, containing $x_\la$ and whose vertices have coordinate $\la,\la-\la_1,\dots,\la-\la_r$. The same argument shows that for any vertex $x\in \kX$, there exists a unique chamber of the building, let say $C_x$, containing $x$ and at minimal distance from $C_-$. Now we claim that 
\begin{eqnarray} \label{lw} \left|V_{\la+w\la_i}(O)\cap
V_{\la_i}(x_\la)\right|=q^{l_w}, \end{eqnarray} where $l_w=d(C_-,C_w)$.
If $l_w=0$ then the claim follows from \eqref{convexite}. So assume that $l_w>0$. Let $G$ be a gallery in $\kP$ of length $l_w$ starting from $C_-$ and finishing in $C_w$. From the preceding discussion we see that there exists $q^{l_w}$ galleries in the building, starting from $C_-$ and whose radial part coincides with $G$. 
Let $(C_-=C_0,C_1,\dots,C_{l_w})$ be one of them. Let $x$ be the unique vertex of $C_{l_w}$ which does not belong to $C_{l_w-1}$. Let us show that $x\in V_{\la+w\la_i}(O)\cap V_{\la_i}(x_\la)$. 
First observe that the sequence of vertices $x_1,\dots,x_{l_w}=x$ defined by $x_k\in C_k$ and $x_k\notin C_{k-1}$ has for radial part the analogue sequence for $G$. In particular $x_\la$ does not belong to this sequence (so $x_\la \in C_x$) and $x\in V_{\la + w\la_i}(O)$. Thus there exists an apartment containing $x_\la$ and $x$, and an isomorphism between this apartment and $\kP$ sending $x_\la$ and $x$ respectively on $\la$ and $\la +w\la_i$. If one composes this isomorphism with the transformation $w^{-1}t_{-\la}$, we obtain an isomorphism sending $x_\la$ on $0$ and $x$ on $\la_i$, which proves that $x\in V_{\la_i}(x_\la)$. Therefore $x\in V_{\la+w\la_i}(O)\cap V_{\la_i}(x_\la)$. Next observe that two different galleries with the same radial part starting from the same chamber cannot have the same last chamber. So they necessarily finish by two different vertices $x$ and $x'$, which proves already one inequality in \eqref{lw}. To show the other one, it suffices to observe that for $x$ in $V_{\la+w\la_i}(O)\cap V_{\la_i}(x_\la)$, the associated chamber $C_x$ is well the last chamber of a gallery starting from $C_-$ with radial part $G$. To see this, notice that \eqref{convexite} implies that any apartment containing $x_\la$ and $O$ contains $C_-$. In particular from the definition of the building we conclude that there exists an apartment containing $C_x$, $C_-$ and $O$. Any copy of $G$ in this apartment is well a gallery starting from $C_-$ finishing in $C_x$ and with radial part $G$. This concludes the proof of \eqref{lw}. Remember that $l(w)$ denotes the length of $w$. It is well known (see \cite{Bou} Ch.6 $\S$ 1.6) that it is also equal to $d(C_+,C_w)$, and that
$$l(w)=\left|\kR^+\cap w^{-1}\kR^-\right|.$$
Moreover, from Theorem 2.15(iv) in \cite{Ron} we get,
$$d(C_-,C_w)+d(C_w,C_+)= d(C_-,C_+)=\left|\kR^+\right|.$$
Thus
$$l_w=|\kR^+|-l(w)=\left|\kR^+\cap w^{-1}\kR^+\right|.$$
Furthermore, by our choice of $w$,  
$$\left|\kR^+\cap w^{-1}\kR^+\right|=\sum_{\alpha \in \kR^+\cap
w^{-1}\kR^+}<\alpha,\la_i>.$$ 
Indeed, for all $\alpha \in
\kR^+$, $<\alpha,\la_i>$ is equal to $0$ or $1$. If there exists 
$\alpha \in \kR^+\cap w^{-1}\kR^+$ such that $<\alpha,\la_i>=0$, then
by decomposing $\alpha$ in the basis $\Pi$, we see that we can always assume $\alpha \in \Pi$. 
But $r_\alpha
\la_i=\la_i$. Thus $w r_\alpha \la_i =w\la_i$. Since $\alpha
\in \Pi$, $r_\alpha$ establishes a bijection from $\kR^+\cap
w^{-1}\kR^+\smallsetminus \{\alpha\}$ onto $\kR^+\cap (w r_\alpha
)^{-1}\kR^+$. Therefore $l(w r_\alpha )= l(w)+1$ and we get an absurdity. We conclude that $<\alpha,\la_i>=1$ for all $\alpha \in \kR^+\cap w^{-1}\kR^+$. 
Since $\alpha \to - \alpha$ is a bijection from
$\kR^+\smallsetminus \kR^+\cap w^{-1}\kR^+$ onto
$w^{-1}\kR^+\smallsetminus \kR^+\cap w^{-1}\kR^+$, we get
$$\sum_{\alpha \in \kR^+\cap
w^{-1}\kR^+}<\alpha,\la_i>=\frac{1}{2}\left\{\sum_{\alpha \in
\kR^+}<\alpha,\la_i>+\sum_{\alpha \in
\kR^+}<\alpha,w\la_i>\right\}.$$ This concludes the proof of the proposition.   
\end{proof}

We can deduce from this proposition the transition probabilities of a radial random walk: 
$$\overline{p}(\la,\la+w\la_i)=q_{t_{\la_i}}^{\frac{1}{2}}\tilde{q}_{t_{w\la_i}}^{\frac{1}{2}}p_i,$$
for all $\la \in P^{++}$. 
In particular for the simple random walk, Formula \eqref{probatransition} gives 
\begin{eqnarray}
\label{trSRW}
\overline{p}(\la,\la+w\la_i)=\frac{\tilde{q}_{t_{w\la_i}}^{\frac{1}{2}}}{\sum_{i=1}^rq_{t_{\la_i}}^{-\frac{1}{2}}N_{\la_i}}.
\end{eqnarray}

\section{Convergence toward the intrinsic BM starting from the interior of the Weyl chamber}      
\label{secdepartbiaise}
Let $x \in
\kX$ and let $(X^x_n,n\ge 0)$ be a nearest neighbor symmetric random walk starting from $x$. Denote by $(Y^x_n,n\ge 0)$ the
$F_0$-random walk starting from $x$.

If $a\in \ka_+$, we denote by $[a]$ one of the weight (it does not matter which one) in $P^+$ at minimal distance from $a$.  
For $N\in \N$ and $t\in \R^+$, we set
$$Y^{N,a}_t= \frac{\overline{Y^{[\sqrt{N}a]}_{[Nt]}}}{\sqrt{N}}.$$
We have 
\begin{theo}
\label{theogeneral} When $N\to +\infty$, the sequence of processes 
$(Y^{N,a}_t,t\ge 0)$ converges in law in 
$\D(\R^+,\overline{\ka_+})$ toward the intrinsic Brownian motion starting from $a$, up to some linear time change. 
\end{theo}
\begin{proof} We first recall that the intrinsic Brownian motion starting from $a\in \ka_+$ takes values a.s. in $\ka_+$. Moreover it is a Feller process with generator $D$ and core $C_c^\infty(\ka_+)$, the space of $C^\infty$ functions with compact support inside $\ka_+$ (see \cite{S} for instance).   
Let $(\kN_t, t\ge 0)$ be a Poisson process with parameter $1$, independent of $(Y_t,t\ge
0)$. Let $(Z^{N,a}_t,t\ge 0)$ be the homogeneous Markov process defined by
$$Z^{N,a}_t= \frac{\overline{Y^{[\sqrt{N}a]}_{\kN_{Nt}}}}{\sqrt{N}} \quad \forall t\ge 0.$$
>From the calculus of the preceding section, we see that the generator $A_N$ of this process is defined for $f\in
C_c^\infty(\ka_+)$ by
$$A_Nf(\la)= \tilde{\rho}^{-1}N\left\{\sum_{i=1}^r q_{t_{\la_i}}^{\frac{1}{2}}p_i
\sum_{\la \in W_0\la_i}
\tilde{q}_{t_{w\la_i}}^{\frac{1}{2}}\frac{F_0(\sqrt{N}\la +
w\la_i)}{F_0(\sqrt{N}\la)}f(\la +
\frac{w\la_i}{\sqrt{N}})\right\}-N f(\la),$$ for all $\la \in
P^+/\sqrt{N}$. Thus \eqref{estimationF02}, implies that for any function $f\in C_c^\infty(\ka_+)$, $A_N f$ converges
uniformly on $\ka_+$ toward $\pi^{-1}p(\partial)(\pi f)$, where 
$p$ is the polynomial 
$$p(x)=\frac{1}{2}\tilde{\rho}^{-1}\left\{\sum_{i=1}^r q_{t_{\la_i}}^{\frac{1}{2}}p_i
\sum_{\la \in W_0\la_i}<\la,x>^2\right\}.$$
In particular $p$ is homogeneous $W_0$-invariant and of degree $2$. Thus it is proportional to $||\cdot||^2$. 
Since $\Delta \pi=0$, $\Delta (\pi f) / \pi = 2Df$. The claim of the theorem follows from Corollary $8.7$ p.$232$ in \cite{EK}.
\end{proof}

\section{Convergence toward the intrinsic BM starting from the origin}

Now we want to obtain a result of convergence toward the intrinsic Brownian motion starting from $0$. It happens to be more difficult to prove than Theorem \ref{theogeneral}, since in this case the criteria of Ethier and Kurtz does not apply directly. Indeed a core of the intrinsic BM starting from $0$ may contain functions which do not vanish in $0$. But for such function $f$ it is not so clear if $A_N f$ converges uniformly toward $Df$ on $\overline{\ka_+}$. However, since the intrinsic BM takes values in $\ka_+$ for $t>0$ (see \cite{ABJ}), the criteria applies for $t\ge \eta$ for any $\eta >0$. Then it only remains to obtain a control of what happens in small time. One way for this is to get good estimates of the transition kernel. But, except when $r=2$, the only known convenient estimates \cite{AST} concern the simple random walk on buildings of type $\tilde{A}_r$. This explains why in the next result we restrict us to this particular case.

We set $Y^N_t:=Y^{N,0}_t$ for all $t\ge 0$, and we denote by $(I_t,t\ge 0)$ the intrinsic BM starting from $0$.
\begin{theo}
\label{theoprincipal} The sequence $(Y^N_t,t\ge 0)$ converges in law
in $\D(\R^+,\overline{\ka_+})$ toward $(I_t,t\ge 0)$, when
$N\to +\infty$.
\end{theo}
\begin{proof} The proof is divided into two parts. We first prove the convergence in law at a fixed time (Lemma \ref{convergenceunid} below). With the criteria of Ethier and Kurtz this allows   
to deduce the convergence in law on $[\eta, +\infty)$ for any $\eta >0$. Then we prove a tightness result, which allows to conclude by general theorems. 

\begin{lem} \label{convergenceunid} For all $t\ge 0$, $Y^N_t$
converges in law toward $I_t$, when $N\to \infty$.
\end{lem}
\begin{proof} Let $f=1_K$ be the indicator function of a compact $K \subset \ka_+$, and let $q_n(O,\cdot)$ be the transition kernel of $(Y_n,n\ge 0)$. Then
\begin{eqnarray*}
\E[f(Y^N_t)]&=&\sum_{\la\in P^+}q_{[Nt]}(O,\la)N_\la f(\frac{\la}{\sqrt{N}})\\
&=& \frac{1}{N^{\frac{r}{2}}}\sum_{u\in \frac{P^+}{\sqrt{N}}\cap
K}N^{\frac{r}{2}}q_{[Nt]}(O,\sqrt{N}u)N_{\sqrt{N}u}.
\end{eqnarray*}
We will show that for all $u\in \frac{P^{++}}{\sqrt{N}}$,
\begin{eqnarray}
\label{convqn}
N^{\frac{r}{2}}q_{[Nt]}(O,\sqrt{N}u)N_{\sqrt{N}u}\to
\text{const}\cdot \frac{\pi(u)^2}{t^{|\kR^+|+\frac{r}{2}}}
e^{-\frac{|u|^2}{2ct}}, \end{eqnarray} when $N\to \infty$. The limit being precisely the density of $I_t$, we will deduce the lemma by a usual argument of integral approximation by series, and Sheff\'e's lemma. First remember that 
$$q_n(O,\la)=p_n(O,\la)\frac{F_0(\la)}{F_0(0)}\tilde{\rho}^{-n},$$
for all $\la\in P^+$ and all $n\ge 0$. Moreover 
by \eqref{formuleintegtransition} and $W_0$-invariance of $h$, we get
$$p_n(O,\la)=\text{const}\cdot
\tilde{\rho}^n q_{t_\la}^{-\frac{1}{2}}\int_U\tilde{h}^n(i\teta)e^{i<\la,\teta>}\ \frac{d\teta}{\mathbf{c}(-i\teta)}.$$
Next we make the change of variables $\teta \to
\teta/\sqrt{[N t]}$ in the integral. Then as in Lemma \ref{convergencepont} we use the convergence of $\tilde{h}^{[N t]}(i\teta /\sqrt{[N t]})$ toward $e^{-\frac{1}{2}||\teta||^2}$, which is moreover dominated by $e^{-\epsilon |\teta|^2}$ for $\epsilon$ small enough and $|\teta| \le \text{const}\cdot \sqrt{[N t]}$. Thus we can apply the dominated convergence theorem, and using 
\eqref{estimationF02}, we get
$$N^{\frac{r}{2}}q_{[Nt]}(O,\sqrt{N}u)N_{\sqrt{N}u} \to \text{const}\cdot \frac{\pi(u)}{t^{\frac{|\kR^+|}{2}+\frac{r}{2}}}
\int
e^{-\frac{1}{2}||\teta||^2+i<\teta,\frac{u}{\sqrt{t}}>}\pi(-i\teta)\ d\teta,
$$
Then we make the change of variables $\teta \to \teta + iu/(c\sqrt{t})$ and we get the limit
$$\text{const}\cdot
\frac{\pi(u)}{t^{\frac{|\kR^+|}{2}+\frac{r}{2}}}e^{-\frac{|u|^2}{2ct}}
\int
e^{-\frac{1}{2}||\teta||^2}\pi(-i\teta+\frac{u}{c\sqrt{t}})\ d\teta.$$
We conclude as follows: observe that the polynomial defined for $x\in \ka$
by 
$$r(x):= \int
e^{-\frac{1}{2}||\teta||^2}\pi(-i\teta+x)\ d\teta.$$ is skew $W_0$-invariant, and of degree exactly $|\kR^+|$. Thus it is equal (up to a constant)
to $\pi$ (see Corollary 3.8
p.362 in \cite{Hel}). This proves \eqref{convqn} and finishes the proof of the Lemma. 
\end{proof}

The next step is the following tightness result. 
\begin{prop}
\label{prop}
Let $\alpha>0$ and $\epsilon>0$. There exists $\eta>0$ such that
$$\limsup_{N\to \infty} \pp\left[\sup_{t\le \eta} |Y^N_t|\ge \alpha\right]\le \epsilon.$$
\end{prop}
\begin{proof} The set 
$$A:=\left\{\sup_{k\le [N\eta]}|Y_k|> \alpha \sqrt{N}\right\},$$
is equal to the disjoint union:
$$A=\coprod_{k=1}^{[N\eta]} A_k,$$
where
$$A_k=\left\{ |Y_1|\le \alpha \sqrt{N}, \dots, |Y_{k-1}|\le \alpha
\sqrt{N},|Y_k|> \alpha \sqrt{N}\right\}.$$ 
For $\epsilon\in(0,1)$, $k\le [N\eta]$, $K>0$
and $\eta\in (0,1)$ such that $\alpha/\sqrt{\eta}>K$, set
$$B_k=\left\{|Y_{[N\eta]-k}|\le
(\frac{\alpha}{\sqrt{\eta}}-K)\sqrt{N\eta} \right\},$$
and
$$D_{\epsilon,k}= \left\{\forall i\in \{([N^{1-\epsilon}]-k)^+,([2N^{1-\epsilon}]-k)^+,\dots,(N-k)^+\},\ \Pi(\overline{Y_i})\ge N^{(|\kR^+|-\epsilon)/2} \right\},$$
where $u^+=u\vee 0$ for any $u\in \R$, and 
$$\Pi(\la):=\prod_{\alpha\in \kR^+} (1+<\alpha,\la>)\quad \forall \la \in P^+.$$
Set also $N_\epsilon:=[N^{1-\epsilon/2|\kR^+|}]$, and 
$$D_\epsilon=\left\{\forall i\in \{[N^{1-\epsilon}],[2N^{1-\epsilon}],\dots,N\}\cap [N_\epsilon,+\infty),\ \Pi(\overline{Y_i})\ge N^{(|\kR^+|-\epsilon)/2} \right\}.$$
Note first that for any $\epsilon\in (0,1)$,
$$\pp[A]\le \pp[B_0^c\cup D_\epsilon^c] + \sum_{k=1}^{N_\epsilon}\pp[A_k]+\sum_{k=N_\epsilon+1}^{[N\eta]}\pp[A_k\cap B_0\cap D_\epsilon].$$
Then Markov property implies that for all $k\ge N_\epsilon$
$$\pp[A_k\cap B_0\cap D_\epsilon] \le \E\left\{1_{A_k}\pp_{Y_k}\left[B_k\cap D_{\epsilon,k}\right]\right\},$$
where for any $x$, $\pp_x$ denotes the law of $Y$ starting from $x$. So 
\begin{eqnarray}
\label{inegalite}
\pp[A]\le \pp[B_0^c\cup D_\epsilon^c] + \left(\sum_{k=1}^{N_\epsilon}\pp[A_k]\right)+\pp[A]\times \sup_{k,x} \pp_x\left[B_k\cap D_{\epsilon,k}\right],
\end{eqnarray}
where the $\sup$ is over integers $k\in(N_\epsilon,[N\eta])$ and vertices $x$ such that $|x|= [\alpha \sqrt{N}]+1$. Now we need the following lemma:  

\begin{lem} 
\label{lemme}
For $\epsilon$ small enough, 
$$\limsup_{N\to +\infty}\ \sup_{k,x}\pp_x\left[B_k\cap D_{\epsilon,k}\right]< 1,$$
where the $\sup$ is over integers $k\le [N\eta]$ and $x$'s such that $|x|=[\alpha\sqrt{N}]+1$.  
\end{lem}
\begin{proof} Let $x$ be such that $|x|=[\alpha\sqrt{N}]+1$. We first claim that
\begin{eqnarray}
\label{eve}
\begin{array}{c}
\text{the probability of the event} \\
\{Y \text{ travel a distance larger than $\alpha\sqrt{N}/2$ in less than $N^{1-\epsilon}$ steps}\}\\
\text{decays exponentially fast to $0$ when $N\to +\infty$}. 
\end{array}
\end{eqnarray}
This is due to the Gaussian factor in the transition kernel of $Y$. More precisely estimates from \cite{AST} show that for any $\la \in P^+$, any $y\in V_\la(x)$ and any $l\ge 1$,   
\begin{eqnarray}
\label{ql}
q_l(x,y)\le C \frac{F_0(\overline{y})F_0(\la)}{F_0(\overline{x})l^{|\kR^+|+r/2}}e^{-c|\la|^2/l}.
\end{eqnarray}
On the other hand for any $\mu\in P^+$, we have 
\begin{eqnarray}
\label{spheres}
|V_\la(x)\cap V_\mu(0)|\le CF_0(\overline{x})q_{t_\la}^{1/2}q_{t_\mu}^{1/2}.
\end{eqnarray}
This follows from \cite{Pa1} Formula (1.7) and Lemma 6.1. 
Then the estimate \eqref{estimationF01} of $F_0$, \eqref{ql} and \eqref{spheres} prove our claim \eqref{eve}.

Now on $D_{\epsilon,k}$ we know that $Y$ (starting from $x$) will enter the set 
$$\kE_{\epsilon,N}:=\{\Pi(\la)\ge N^{(|\kR^+|-\epsilon)/2}\}$$ 
before time $N^{1-\epsilon}$. 
So if we denote by $x'$ the position of $Y$ at his first entrance in this set, we can assume by \eqref{eve} that $|x'|\ge \alpha\sqrt{N}/2$. By using next the Markov property, it is enough to prove the lemma with $x'$ in place of $x$. Actually in the following, for notation convenience, we will assume that $x$ already lies in  the set $\kE_{\epsilon,N}$.

Then let $T_{\partial \ka_+}$ be the first time the (radial part of) the random walk reaches $\partial \ka_+$. We claim that for $\epsilon$ small enough and $N$ large enough,
\begin{eqnarray}
\label{murs}
\pp_x\left[B_k\cap D_{\epsilon,k}\cap\{T_{\partial \ka_+}\le [N\eta]-k\}\right]\le \frac{1}{4} \quad \forall k \le [N\eta].
\end{eqnarray}
To see this we need to describe the behavior of $\overline{Y}$ on $P^{++}$. Denote by $\overline{q}(\cdot,\cdot)$ its one-step transition kernel, and remember that 
$\overline{p}(\cdot,\cdot)$ denotes the one of the radial part of the simple random walk on $\kX$. First by definition
$$\overline{q}(\la,\mu)=\frac{F_0(\mu)}{F_0(\la)}\tilde{\rho}^{-1}\overline{p}(\la,\mu),$$
for all $\la \in P^{++}$ and $\mu\in P^+$ such that $|\la-\mu|=1$. Denote by $c_0$ the inverse of the number of neighbors (in $P^+$) of any weight $\la \in
P^{++}$. By using the explicit formula \eqref{trSRW} of $\overline{p}(\la,\mu)$, we see that 
\begin{eqnarray}
\label{qp}
\overline{q}(\la,\mu)= c_0 \frac{F_0(\mu)q_{t_\mu}^{1/2}}{F_0(\la)q_{t_\la}^{1/2}}.
\end{eqnarray}
Let now $y\in \kX$ be such that $\overline{y}\in \kE_{\epsilon,N}$ and 
$|y|\le N^{(1+\epsilon/|\kR^+|)/2}$. If $\la=\overline{y}$, then by \eqref{qp} we get
\begin{eqnarray}
\label{mula}
\pp_y\left[T_{\partial \ka_+}\le N^{1-\epsilon}\right]=\sum \frac{F_0(\mu)q_{t_\mu}^{1/2}}{F_0(\la)q_{t_\la}^{1/2}}c_0^l,
\end{eqnarray}
where the sum is over all $l\le N^{1-\epsilon}$ and all paths $(\la=\la_1,\dots,\la_l=\mu)$ in $P^+$ of length $l$, going from $\la$ to some $\mu\in \partial\ka_+\cap P^+$. But if $\mu \in \partial \ka_+\cap P^+$, then at least one positive root vanishes in $\mu$. If moreover $\mu$ is at distance at most $N^{1-\epsilon}$ from $\la=\overline{y}$, then $|\mu|\le C |\la|$ and thus
$$
\Pi(\mu)\le C N^{(|\kR^+|+\epsilon-1)/2}.
$$
So in \eqref{mula} the ration $F_0(\mu)q_{t_\mu}^{1/2}/F_0(\la)q_{t_\la}^{1/2}$ is always bounded by $CN^{\epsilon-1/2}$. 
This shows that 
\begin{eqnarray}
\label{T}
\pp_y\left[T_{\partial \ka_+}\le N^{1-\epsilon}\right] &\le & \nonumber C N^{\epsilon-1/2}\overline{P}(T_{\partial \ka_+}\le N^{1-\epsilon})\\
                                                       &\le & C N^{\epsilon-1/2},
\end{eqnarray}                                                
where $\overline{P}$ denotes the law of the simple random walk on $P^+$ (say reflected on $\partial \ka_+$). 
On the other hand, thanks again to the Gaussian factor in \eqref{ql}, we know that (except on an event of probability decaying exponentially fast to $0$ when $N\to +\infty$) for all $i\in \{([N^{1-\epsilon}]-k)^+,([2N^{1-\epsilon}]-k)^+,\dots,(N-k)^+\}$, 
$|Y_i|\le N^{(1+\epsilon/|\kR^+|)/2}$. Thus \eqref{T}, applied (at most) $[\eta N^\epsilon]$ times, proves our claim \eqref{murs}, 
if $\epsilon$ is small enough (e.g. $\epsilon<1/4$).

Consider now $|\mu|\le |x|$, and some path $f=(\overline{x}=\la_1,\dots,\la_l=\mu)$ going from $\overline{x}$ to $\mu$ and not intersecting 
$\partial \ka_+$. Set $\alpha_0:=\alpha_1+\dots+\alpha_r$, and let $H_x$ be the affine hyperplane containing $x$ and orthogonal to $\alpha_0$. 
Call $f'$ the path obtained by reflecting orthogonally the last excursion of $f$ out of $H_x$ on the other side of $H_x$. Then $f'$ does not intersect $\partial \ka_+$ and the map $f\mapsto f'$ is one to one. 
Moreover under $\overline{P}$, $f$ and $f'$ have the same probability. 
Thus by using again \eqref{qp} and \eqref{estimationF02}, we see that the probability (for $\overline{Y}$) to arrive at some time $t$ below level $|x|/2$ without touching $\partial \ka_+$ is at most $1/2+o(1)$, where $o(1)$ tends to $0$ when $N\to +\infty$. Thus for $N$ large enough,
\begin{eqnarray*}
\pp_x\left[B_k\cap D_{\epsilon,k}\cap\{T_{\partial \ka_+}\ge [N\eta]-k\}\right]\le \frac{2}{3} \quad \forall k \le [N\eta].
\end{eqnarray*}  
This concludes the proof of the lemma. 
\end{proof}
 
We can finish now the proof of Proposition \ref{prop}. Coming back to \eqref{inegalite} and applying the previous lemma, we see that for $\epsilon$ small enough and $N$ large enough, 
$$\pp[A] \le C\left(\pp[B_0^c]+\pp[D_\epsilon^c]+\sum_{k=1}^{N_\epsilon}\pp[A_k]\right).$$
But the cardinality of a sphere of radius $\la$ is of order $q_{t_\la}$. So estimate \eqref{ql} with $x=O$ (and $\overline{y}=\la$) shows that the right hand term in the above inequality tends to $0$ when $K$ and $N$ tend to $+\infty$. This finishes the proof of the proposition. \end{proof} 

\noindent Theorem \ref{theoprincipal} is now a consequence of standard theorems (see for instance \cite{Bil} or Theorem $7.2$ p.$128$ in \cite{EK}).
\end{proof}

\begin{rem} \emph{As mentioned previously, when $r=2$ the estimates in \cite{AST} hold for all symmetric nearest neighbor random walks. So slight modifications of the previous proof show that Theorem \ref{theoprincipal} extend to all these random walks when $r=2$.}
\end{rem}

\textit{This paper is part of my Ph.D. thesis. I warmly thank my advisors Jean-Philippe Anker and Philippe Bougerol for their suggestion of studying this problem, and Philippe Bougerol for having drawn my attention to a mistake in the proof of Proposition \ref{prop} in a previous version of this paper.}

\vspace{0.3cm}
\noindent \textit{Universit\'e d'Orl\'eans, F\'ed\'eration Denis Poisson, Laboratoire MAPMO \\
B.P. 6759, 45067 Orl\'eans cedex 2, France.}\\
\noindent mail: bruno.schapira@math.u-psud.fr

\vspace{0.2cm}
\noindent \textit{Current address: D\'epartement de Math\'ematiques, B\^at. 425, Universit\'e Paris-Sud, F-91405 Orsay
cedex, France. }

\end{document}